\documentclass[english,reqno]{amsart}
\usepackage{etex}
\usepackage{amsmath,amssymb,amsthm,bbm,mathtools,comment}
\usepackage[shortlabels]{enumitem}
\usepackage[pdftex,colorlinks,backref=page,citecolor=blue]{hyperref}
\hypersetup{pdfpagemode=UseNone,pdfstartview={XYZ null null 1.00}}
\usepackage[mathscr]{euscript}
\usepackage[usenames,dvipsnames]{color}
\usepackage{adjustbox,tikz,calc,graphics,babel,standalone}
\usetikzlibrary{shapes.misc,calc,intersections,patterns,decorations.pathreplacing}
\usetikzlibrary{arrows,shapes,positioning}
\usetikzlibrary{decorations.markings}
\usepackage[final]{microtype}
\usepackage[numbers]{natbib}
\usepackage{cmtiup}
\usepackage{amsfonts}
\usepackage{graphicx}
\usepackage{caption}
\usepackage{subcaption}
\usepackage{verbatim}
\usepackage{array}
\usepackage[frame,cmtip,arrow,matrix,line,graph,curve]{xy}
\usepackage{graphpap, color, pstricks}
\usepackage{pifont}
\usepackage[final]{microtype}
\usepackage{cmtiup}

\allowdisplaybreaks

\theoremstyle{plain}
\newtheorem{theorem}{Theorem}[section]		
\newtheorem{lemma}[theorem]{Lemma}
\newtheorem{claim}[theorem]{Claim}

\theoremstyle{remark}

\def\CC{\mathscr{C}}

\def\bd{7k}

\def\N{\mathbb{N}}

\let\originalleft\left
\let\originalright\right
\renewcommand{\left}{\mathopen{}\mathclose\bgroup\originalleft}
\renewcommand{\right}{\aftergroup\egroup\originalright}

\makeatletter
\def\imod#1{\allowbreak\mkern10mu({\operator@font mod}\,\,#1)}
\makeatother

\begin{document}
\title{Subgraphs of large connectivity and chromatic number}

\author{Ant\'onio Gir\~ao}
\address{School of Mathematics, University of Birmingham, Edgbaston, Birmingham, B15\thinspace2TT, UK}
\email{giraoa@bham.ac.uk}

\author{Bhargav Narayanan}
\address{Department of Mathematics, Rutgers University, Piscataway, NJ 08854, USA}
\email{narayanan@math.rutgers.edu}

\date{2 April 2020}

\begin{abstract}
Resolving a problem raised by Norin, we show that for each $k \in \N$, there exists an $f(k) \le \bd$ such that every graph $G$ with chromatic number at least $f(k)+1$ contains a subgraph $H$ with both connectivity and chromatic number at least $k$. This result is best-possible up to multiplicative constants, and sharpens earlier results of Alon--Kleitman--Thomassen--Saks--Seymour from 1987 showing that $f(k) = O(k^3)$, and of Chudnovsky--Penev--Scott--Trotignon from 2013 showing that $f(k) = O(k^2)$. Our methods are robust enough to handle list colouring as well: we also show that for each $k \in \N$, there exists an $f_\ell(k) \le 4k$ such that every graph $G$ with list chromatic number at least $f_\ell(k)+1$ contains a subgraph $H$ with both connectivity and list chromatic number at least $k$. This result is again best-possible up to multiplicative constants; here, unlike with $f(\cdot)$, even the existence of $f_\ell(\cdot)$ appears to have been previously unknown.
\end{abstract}

\maketitle

\section{Introduction}
Many of the central open problems in graph theory concern structures that are unavoidable in graphs of large chromatic number, Hadwiger's conjecture~\cite{Hadwiger} being perhaps the most notable example. 

Here, we shall be concerned with an extension of the following well-known fact: every graph of chromatic number at least $4k+1$ contains a subgraph of connectivity at least $k$, as follows from a classical result of Mader~\citep{Mader} asserting that every graph of minimum degree at least $4k$ contains a $k$-connected subgraph. It is natural to then ask if a graph of large chromatic number must contain a subgraph of both large connectivity and large chromatic number; this was answered by Alon, Kleitman, Thomassen, Saks and Seymour~\citep{Alon} who showed for each $k \in \N$ that there exits an $f(k) = O(k^3)$ such that every graph with chromatic number at least $f(k)+1$ contains a subgraph whose connectivity and chromatic number are both at least $k$. This was improved by Chudnovsky, Penev, Scott and Trotignon~\citep{Alex} who (amongst other things) showed that $f(k) = O(k^2)$, and the lower order terms in this result were later improved by Penev, Thomass\'e and Trotignon~\citep{Thomasse}.

The results described above have since found many applications in the study of graphs of large chromatic number. Motivated by applications to the study  of Hadwiger's conjecture (see~\citep{H1, H2}), Norin~\citep{ques} asked if the aforementioned results could be sharpened to show an essentially best-possible estimate of $f(k) = O(k)$; our first result answers this question affirmatively.

\begin{theorem}\label{mainthm}
For each $k \in \N$, every graph $G$ with chromatic number at least $\bd+1$ contains a subgraph $H$ with both connectivity and chromatic number at least $k$.
\end{theorem}

In other words, Theorem~\ref{mainthm} asserts that $f(k) \le 7k$, and from below, Alon, Kleitman, Thomassen, Saks and Seymour~\citep{Alon} showed that $f(k) \ge 2k-3$. While these bounds are not too far apart, we make no particular effort to optimise the multiplicative constant in our result since it seems unlikely that this will completely bridge the gap between the upper and lower bounds.

With the analogue of Hadwiger's conjecture for list colourings~\citep{lists} in mind, it is also natural to ask if, for each $k \in \N$, there exits an $f_\ell(k) \in \N$ such that every graph with list chromatic number at least $f_\ell(k)+1$ contains a subgraph whose connectivity and list chromatic number are both at least $k$. The methods we use to prove Theorem~\ref{mainthm} are robust enough to answer this question; our second result, stated below, is again essentially best-possible.

\begin{theorem}\label{mainthm2}
For each $k \in \N$, every graph $G$ with list chromatic number at least $4k+1$ contains a subgraph $H$ with both connectivity and list chromatic number at least $k$.
\end{theorem}
In other words, Theorem~\ref{mainthm2} asserts that $f_\ell(k) \le 4k$, and the construction of Alon, Kleitman, Thomassen, Saks and Seymour~\citep{Alon} mentioned earlier also shows that $f_\ell(k) \ge 2k-3$.

This paper is organised as follows. We introduce the key notions that we need in Section~\ref{s:prelim}, give the proofs of both our results in Section~\ref{s:proof}, and conclude with a discussion of some related questions in Section~\ref{s:conc}.

\section{Preliminaries}\label{s:prelim}
We start by establishing some notation. For a set $X$, we write $2^X$ for the power set of $X$, and given a function $\lambda$ defined on $X$ and a subset $Y \subset X$, we write $\lambda \, |_Y$ for the restriction of $\lambda$ to $Y$. Given a graph $G$, as is usual, we write $\chi(G)$, $\chi_\ell(G)$ and $\kappa(G)$ for the chromatic number, list chromatic number and the connectivity of $G$ respectively, and for a subset $X \subset V(G)$, we write $G[X]$ for the subgraph of $G$ induced by $X$; for any graph-theoretic terminology not defined here, we refer the reader to~\citep{book}.

In what follows, we fix $k \in \N$ and work with a fixed palette $\CC$ of $\bd$ colours. Our proof of Theorem~\ref{mainthm} will hinge around two notions, those of templates and extensibility, that we define below; the modifications needed for Theorem~\ref{mainthm2} will be indicated in place where needed, in Section~\ref{s:proof}.

\subsection*{Templates}
A template on a graph consists of a set of properly `pre-coloured' vertices, along with lists of `forbidden' colours at each of the remaining vertices. Formally, a \emph{template $T = (S, c, F)$} on a graph $G$ consists of 
\begin{enumerate}
\item a subset $S \subset V(G)$, 
\item a proper \emph{pre-colouring} $c:S \to \CC$ of the induced subgraph $G[S]$, and 
\item a function $F: V(G)\setminus S \to 2^\CC$ specifying a list of \emph{forbidden} colours at each remaining vertex.
\end{enumerate}

Let $G = (V, E)$ be a graph and let $T = (S, c, F)$ be a template on $G$. We define the \emph{degree} of $T$ by
\[ \deg(T) = k|S| + \sum_{v \in V \setminus S} |F(v)|.\]
For a set $X \subset V$, the \emph{restriction} of $T$ to the induced subgraph $G[X]$ is naturally the template 
\[T_X = (S \cap X, c\,|_{S \cap X}, F\,|_{S \cap X}).\] 
Let us note that the degree of a template is additive across disjoint restrictions, i.e., if $X \cup Y$ is a partition of $V$, then 
\[\deg(T) = \deg(T_X) + \deg(T_Y).\] 
Finally, we say that a proper colouring $\hat c:V \to \CC$ of $G$ \emph{respects $T$} if it extends the pre-colouring specified by $T$ while avoiding the forbidden colours at all the other vertices, i.e., if
\begin{enumerate}
\item $\hat c(v) = c(v)$ for all $ v \in S$, and
\item $\hat c(v) \notin F (v)$ for all $v \in V \setminus S$.
\end{enumerate}

\subsection*{Extensibility}
We say that a graph $G$ is \emph{inextensible} if there exists a template $T = (S, c, F)$ on $G$ such that
\begin{enumerate}
\item $\deg(T) \le 2k^2$,
\item $|F(v)| \le 2k$ for each $v \in V(G) \setminus S$, and
\item there is no proper colouring of $G$ using the palette $\CC$ that respects $T$;
\end{enumerate}
we call this template $T$ a \emph{witness} for the inextensibility of $G$, and also note that there may be multiple templates satisfying the requisite conditions and witnessing this fact. Analogously, we say that $G$ is \emph{extensible} if it is not inextensible. 

First, we observe that graphs of sufficiently large chromatic number are inextensible.
\begin{lemma}\label{inext}
If $G$ is a graph with $\chi(G) \ge \bd + 1$, then $G$ is inextensible.
\end{lemma}
\begin{proof}
This is obvious; the empty template, with no colours forbidden at any vertex and no pre-coloured vertices, shows that $G$ is inextensible.
\end{proof}

Next, there is some slack in the definition of inextensibility, as shown by the following observation.

\begin{lemma}\label{strongwitness}
If a graph $G$ is inextensible, then this is witnessed by a template $T = (S, c, F)$ on $G$ for which $|F(v)| \le k-1$ for each $v \in V(G) \setminus S$.
\end{lemma}
\begin{proof}
Of all the templates $T = (S, c, F)$ witnessing the inextensibility of $G$, choose one with $|S|$ maximal, and suppose for the sake of contradiction that there exists a vertex $v \in V(G) \setminus S$ with $|F(v)| \ge k$. We claim that we may find another template $T'$ witnessing the inextensibility of $G$ in which $S \cup \{x\}$ is pre-coloured, contradicting the maximality of $|S|$. 

To see this, first note that since $\deg(T) \le 2k^2$, we must have $|S| \le 2k$. Next, as $|F(v)| \le 2k$ and $|\CC| = \bd$, there exists a colour in $\CC \setminus F(v)$ not appearing anywhere in the pre-colouring of $S$; to obtain $T'$ from $T$, we pre-colour $v$ with this colour and remove the forbidden list of colours associated with $v$; any proper colouring of $G$ respecting $T'$ also respects $T$, so it suffices to check that $\deg(T') \le \deg(T) \le 2k^2$. This is straightforward: in passing from $T$ to $T'$, the first summand in the definition of the degree increases by $k$, but since $|F(v)| \ge k$, the second summand decreases by at least $k$.
\end{proof}

Of course, we could alter the definition of inextensibility to remove this slack, but we prefer the definition above since the elbow room makes our inductive argument more transparent.

\section{Proofs of the main results}\label{s:proof}

With the notions of templates and extensibility in hand, we are now ready to prove our first main result.

\begin{proof}[Proof of Theorem~\ref{mainthm}]

Let $G$ be a graph with $\chi(G) \ge \bd + 1$ which, by Lemma~\ref{inext}, is inextensible. Our goal is to find a subgraph $H$ of $G$ with both $\chi(H) \ge k$ and $\kappa (H) \ge k$. 

Let $H$ be a \emph{minimal inextensible} induced subgraph of $G$ on some vertex set $U$, and let $T = (S,c,F)$ be a template on $H$ witnessing its inextensibility with $|F(v)| \le k-1$ for each $v \in U \setminus S$, as promised by Lemma~\ref{strongwitness}. As observed earlier, since $\deg(T) \le 2k^2$, it must be the case that $|S| \le 2k$.

We first show that $H$ has large connectivity.

\begin{claim}\label{c1}
$\kappa(H) \ge k$.
\end{claim}
\begin{proof}
Suppose for the sake of contradiction that $H$ is not $k$-connected. We shall find a proper colouring $\hat c$ of $H$ using  $\CC$ that respects $T$, contradicting the inextensibility of $H$ as witnessed by $T$.

First, if $H$ is isomorphic to a complete graph $K_k$, then the construction of $\hat c$ is straightforward. It suffices to find a proper colouring of $H[U \setminus S]$ where each vertex $v$ receives a colour from the list $L_v$ of colours in $\CC$ appearing neither somewhere in $S$, nor in $F(v)$. Since $|F(v)| \le k-1$, we see that for each vertex $v \in U \setminus S$, we have (with room to spare)
\[|L_v| \ge |\CC| - |S| - |F(v)| \ge \bd -2k - (k-1) \ge k.\]
Since $H$ has $k$ vertices and $|L_v| \ge k$ for each $v \in U \setminus S$, we may colour each vertex $v \in U \setminus S$ with a colour from $L_v$ in such a way that these vertices all get distinct colours.

Next, suppose that there is a subset $X \subset U$ of size at most $k-1$ which disconnects $H$, and fix a nontrivial partition $Y \cup Z$ of $
U \setminus X$ where there are no edges between $Y$ and $Z$. Since
\[ \deg(T_Y) + \deg(T_Z) \le \deg(T) \le 2k^2,\]
we assume, without loss of generality, that $\deg(T_Z) \le k^2$.

First, let $U' = X \cup Y$ and consider $H' = H[U']$. Starting with the template $T_{U'}$, we construct a new template $T' = (S \cap U', c\,|_{S \cap U'}, F')$ on $H'$ by defining $F'$ as follows: we start with $F'(v) = F(v)$ for each $v \in U'$, and then for each $z \in S \cap Z$, we add the colour $c(z)$ to $F'(v)$ for each neighbour $v \in X$ of $z$. 

It is easy to see that $|F'(v)| \le 2k$ for each $v \in U'$; indeed, 
\[|F'(y)| = |F(y)| \le k-1\] 
for each $y \in Y$, and since $\deg(T_Z) \le k^2$, we must have $|S \cap Z| \le k$, so
\[|F'(x)| \le |S \cap Z|+ |F(x)| \le k + (k-1) = 2k -1\] 
for each $x \in X$. 

It is also not hard to see that $\deg(T') \le \deg(T)$; indeed, in passing from $T$ to $T'$, the removal of the pre-coloured vertices in $S\cap Z$ decreases the degree of the template by $k|S \cap Z|$, while the addition of the colours of these vertices to the lists of forbidden colours in $X$ increases the degree of the template by at most $|X||S \cap Z| \le (k-1)|S \cap Z|$.

From the minimality of $H$, we know that $H'$ is extensible, so there exists a proper colouring $c'$ of $H'$ using $\CC$ that respects $T'$. 

Next, we take $U'' = X \cap Z$ and $H'' = H[U'']$, and construct a new template $T''$ on $H''$ as follows: we start with the restriction $T_Z$, and then additionally pre-colour the vertices in $X \setminus S$ according to $c'$. Clearly, we have
\[ \deg(T'') \le k|X| + \deg(T_Z) \le k(k-1) + k^2 \le 2k^2,\]
and since all the lists of forbidden colours in $T''$ are inherited from $T$, these lists all have size at most $k-1$. 

From the minimality of $H$, we again know that $H''$ is extensible, so we may find a proper colouring $c''$ of $H''$ using $\CC$ that respects $T''$. 

Of course, we have ensured that $c'\,|_X = c''\,|_X$, so gluing these two colourings together along $X$ gives us a proper colouring $\hat c$ of $H$ as required.
\end{proof}

Next, we show that $H$ has large chromatic number.

\begin{claim}\label{c2}
$\chi(H) \ge k$.
\end{claim}
\begin{proof}
Suppose again for the sake of contradiction that $\chi(H) \le k-1$. This certainly means that we may partition $U \setminus S$ into $k-1$ independent sets $J_1, J_2, \dots, J_{k-1}$.

Now, we order the vertices of $U \setminus S$ in such a way that each independent set $J_i$ forms an interval in this ordering. We process the vertices of $U \setminus S$ in order and partition them into (at most) $3k$ intervals $I_1, I_2, \dots, I_{3k}$, each contained within some independent set $J_i$, as follows: having constructed $I_1, I_2, \dots, I_{m-1}$, we consider as yet unprocessed vertices in order and add them one by one to $I_m$, and when considering a vertex $v$, we decide to stop the construction of $I_m$ and move on to $I_{m+1}$ based on the following pair of rules:
\begin{enumerate}
\item we stop \emph{without adding} $v$ to $I_m$ if the addition of $v$ stops $I_m$ from being contained within a single independent set $J_i$, and otherwise
\item we stop \emph{by adding} $v$ to $I_m$ if this causes the sum $\sum_{u \in I_m} |F(u)|$ to exceed $k$.
\end{enumerate}

It is not hard to see that this procedure produces at most $3k$ intervals: the number of times we stop on account of the first rule is at most $k-1$, and the number of times we stop on account of the second rule is at most $2k$ since $\deg(T) \le 2k^2$. Notice that these intervals $I_1, I_2, \dots, I_{3k}$ have the following properties:
\begin{enumerate}
\item each $I_m$ is an independent set in $H$, and
\item each $I_m$ satisfies $\sum_{v \in I_m}|F(v)| \le 2k$.
\end{enumerate}

We now construct a proper colouring $\hat c$ of $H$ using $\CC$ that respects $T$, again contradicting the inextensibility of $H$ as witnessed by $T$. For each interval $I_m$, we form the list $L_m$ of colours appearing neither somewhere in $S$, nor in $\cup_{v \in I_m} F(v)$, noting that 
\[|L_m| \ge |\CC| - |S| - \sum_{v \in L_m} |F(v)| \ge \bd -2k - 2k \ge 3k.\] 
Since there are at most $3k$ intervals, it is now clear, as before, that we can each interval $I_m$ with some colour from $L_m$ in such a way that distinct intervals get distinct colours. This colouring yields a proper colouring $\hat c$ of $H$ as required.
\end{proof}

Claims~\ref{c1} and~\ref{c2} together show that $H$ has the requisite properties, proving the result.
\end{proof}

The proof of our second main result is almost identical to, and easier in parts than, the proof of Theorem~\ref{mainthm}, so we only outline the required modifications.
\begin{proof}[Proof of Theorem~\ref{mainthm2}]
The main point here is to tweak the notion of a template to deal with list colourings. This is straightforward: given a graph with lists at each vertex, a template now consists of some vertices properly pre-coloured with colours from their lists, along with some colours forbidden from the lists at each of the remaining vertices. With this modification in place, the rest of the proof is more or less identical to that of Theorem~\ref{mainthm}; the analogue of Claim~\ref{c2} becomes easier, and on account of this, the argument allows us to start with a bound of $\chi_\ell(G) \ge 4k+1$ (in place of the bound of $\chi(G) \ge 7k+1$ that is needed in the proof of Theorem~\ref{mainthm}).
\end{proof}

\section{Conclusion}\label{s:conc}
We conclude with a discussion of some problems that still remain. In what follows, we restrict ourselves to proper colourings, but similar considerations also apply to list colourings. 

As mentioned earlier, we now know that $2k-3 \le f(k) \le 7k$. There are reasons to believe that the lower bound is more reflective of the truth, as we  now explain.

All of~\citep{Alon, Alex, Thomasse} study an asymmetric analogue of the problem treated here: for $m,k \in \N$, let $g(m,k)$ be the least natural number such that every graph $G$ with chromatic number at least $g(m,k)+1$ contains a subgraph $H$ with connectivity at least $k$ and chromatic number at least $m$. 

Of course, it is clear from our results that $g(m,k) = \Theta(m+k)$, but more precise results are available in the `off-diagonal' case when $m$ is much larger than $k$: it is shown in~\citep{Alon} that $g(m,k) \ge m+k-3$, and~\citep{Thomasse} shows that $g(m,k) \le m+2k-3$ when $m \ge 2k^2$. It is not hard to modify the arguments here to establish something like this latter bound unconditionally, i.e., to show that $g(m,k)  = m + O(k)$ for all $m, k \in \N$. Obtaining a precise description of $g(m,k)$ when $m$ is much larger than $k$ may be a good starting point towards pinning down the exact value of $f(k) = g(k,k)$.

\section*{Acknowledgements}
The first author wishes to acknowledge support from EPSRC grant EP/N019504/1, and the second author was partially supported by NSF grant DMS-1800521. We would like to thank Sergey Norin and Sophie Spirkl for several helpful discussions.

\bibliographystyle{amsplain}
\bibliography{kappa_chi}

\end{document}